\documentclass[12pt]{article}

\usepackage{typearea}  
\usepackage{amsmath}  
\usepackage{amsfonts}  
\usepackage[english]{babel}  
\usepackage[T1]{fontenc}   

\newtheorem{prop}{Proposition}
\newtheorem{lemma}{Lemma}
\newtheorem{Rem}{Remark}

\newcommand{\CQFD}{\nolinebreak\hfill\rule{2mm}{2mm}\medbreak\par}  
  
\newcommand{\ind}{\mathbf{1}}  
  
\newcommand{\disp}{\displaystyle}

\newcommand{\s}{\sigma(S^{d-1})}  
\newcommand{\om}{\omega}

\def\rit{\mathbb{R}}  
  
\def\nit{\mathbb{N}}

\def\E{\mathop{\hbox{\rm I\kern-0.20em E}}\nolimits}

\def\med{\hskip 10pt}  
\def\lng{\hskip 20pt}  
  
\def\og{\leavevmode\raise.3ex  
     \hbox{$\scriptscriptstyle\langle\!\langle$~}}  
\def\fg{\leavevmode\raise.3ex  
     \hbox{~$\!\scriptscriptstyle\,\rangle\!\rangle$}~}

\allowdisplaybreaks  
  
\title{On finite range stable-type concentration}  

\author{
J.-C. Breton$^{(1,a)}$\setcounter{footnote}{-1}\footnote{(a) Email: jcbreton@univ-lr.fr, I thank the School of Mathematics at the Georgia Institute of Technology where part of this research was done.  
  }  
 and C. Houdr\'e$^{(2,3,b)}$\setcounter{footnote}{-1}\footnote{(b) Email: houdre@math.gatech.edu, Research supported in part by the NSA grant G-37-A56/C.}\\
\\
\small{(1) Laboratoire de Mathématiques et Applications,}\\
\small{Universit\'e de La Rochelle, 17042 La Rochelle Cedex, France}\\
\\
\small{(2) Laboratoire d'Analyse et de Mathématiques Appliquées, CNRS UMR 8050}\\
\small{Université Paris XII, 94010 Créteil Cedex, France }\\
\\
\small{(3) School of Mathematics,}\\
\small{Georgia Institute of Technology, Atlanta, GA 30332, U.S.A.}\\
}

\begin{document}  
  
\maketitle  
\begin{abstract}   
 We study the deviation probability $P(f(X)-E[f(X)]\geq x)$ where $f$ is a Lipschitz (for the Euclidean norm) function  defined on $\rit^d$ 
 and $X$ is an $\alpha$-stable random vector of index $\alpha \in (1,2)$. We show that the order of this deviation is $e^{-cx^{\alpha/(\alpha-1)}}$ 
 or $e^{-cx^\alpha}$ according to $x$ taking small values or being in a finite range interval.
 In the second part, we generalizes these finite range deviations to $ P(F-m(F)\geq x)$ where $F$ is a stochastic functional on the Poisson space equipped with a   
 stable Lévy measure of index $\alpha\in(0,2)$ and where $m(F)$ is a median of $F$.\\

\vskip 5pt 
\noindent
{\bf AMS 2000 Suject Classification.} 60G70, 62G07, 62C20, 41A25.  
 \end{abstract}


\section{Introduction and preliminaries}  
\label{sec:preliminary}

The purpose of these notes is to further complete our understanding of the stable concentration phenomenon, 
by obtaining the finite range behavior of $P(F-E[F]\geq x)$, with $F=f(X)$ where $f$ is a Lipschitz function and $X$ is a stable random vector
 or with $F$ a stochastic functional on the Poisson space equipped with a stable L\'evy measure. 
 
\vskip 5pt  
  
More precisely, in Section \ref{sec:M}, we consider an $\alpha$-stable random vector $X\in \rit^d$ of index $\alpha \in(1,2)$ and $f:\rit^d\to\rit$ a $1$-Lipschitz function   
  and we investigate deviation probabilities with respect to the mean,
\begin{equation}
\label{eq:dev-def}
P(f(X)-E[f(X)]\geq x),
\end{equation}
for all order of $x$. When $x$ takes small values, the deviation \eqref{eq:dev-def} is bounded by $e^{-cx^{\alpha/(\alpha-1)}}$, while in Section \ref{sec:intermediaire1} and \ref{sec:intermediaire2}, 
we give results for intermediate values of $x$, in which case the bound is of order $e^{-cx^\alpha}$.  In the intermediate range case, we extend the result of \cite{M2} which hold for $\alpha$ close enough to $2$, to any $\alpha\in(1,2)$. These results complement the $(1/x^\alpha)$--behavior given for large values of $x$ in \cite[Th. 1, Th. 2]{HM}   
 and \cite[Th. 6.1]{BHP} and generalize \cite[Th. 2]{M2} and \cite[Th. 6.3]{BHP}. 
  
  In Section \ref{sec:BHP}, we further extend the intermediate range results to stochastic functionals on the Poisson space $\Omega^{\rit^d}$ of $\rit^d$ equipped with a stable Lévy measure.
In this case, the deviations are expressed with respect to a median of the stochastic functionals and we recover as a special case the deviation of Lipschitz
function of stable vector of any index $\alpha\in (0,2)$. We also briefly state the small values result for functionals when $\alpha>1$.

The main ideas of the proof are present in Section \ref{sec:M}. In Section \ref{sec:BHP}, dealing with a median rather than with the mean makes the argument more involved, 
this is the reason for dealing first with stable random vector of index $\alpha\in (1,2)$  in Section \ref{sec:M}.

The general scheme of proof is to decompose the stable random vector $X$ in two independent components $Y_R+Z_R$ with $R>0$ a level of truncation for the L\'evy measure. $R$ is then chosen according to the range of deviation we are interested in.  
We argue similarly in Section \ref{sec:BHP} with truncated configuration $\omega_R=\omega \cap B(0,R)$, $\omega\in\Omega^{\rit^d}$, where $B(0,R)$ is the (Euclidean) ball of radius $R$ centered at $0$. In both cases, the remainder is controlled alternatively.


\section{Stable random vector of index $\alpha>1$}  
\label{sec:M}

  In this part, we consider an $\alpha$-stable random vector $X\in \rit^d$ of index $\alpha \in(1,2)$ and $f:\rit^d\to\rit$ a $1$-Lipschitz function with respect to Euclidean norm $\|\cdot\|$   and we investigate the deviation \eqref{eq:dev-def}. 
The vector $X$ has an infinitely divisible law $ID(b,0,\nu)$ with L\'evy measure given by   
\begin{equation}  
\label{eq:nu}  
\nu(B)=\int_{S^{d-1}}\sigma(d\xi) \int_0^\infty \ind_B(r\xi) r^{-1-\alpha} dr, \quad B\in{\cal B}(\rit^d).  
\end{equation}  
Its characteristic function is $\varphi_X=e^{\phi_X}$ with 
$$
\phi_{X}(u)=i\langle u,b\rangle+\int_{(\rit^*)^d}\big(e^{i\langle u,y\rangle}-1-i\langle u,y\rangle\ind_{||y||\leq 1}\big) \nu(dy),
$$
where $(\rit^*)^d=\rit^d\setminus \{0\}$. Write $X=Y_R+Z_R$ where $Y_R$ and $Z_R$ are two independent infinitely divisible random vectors with respective characteristic functions 
 $\varphi_{Y_R}=e^{\phi_{Y_R}}$ and $\varphi_{Z_R}=e^{\phi_{Z_R}}$ whose exponents are given by 
 \begin{eqnarray*}
 \phi_{Z_R}(u)&=&\int_{||y||>R}\big(e^{i\langle u,y\rangle}-1\big) \nu(dy)\\
\phi_{Y_R}(u)&=&i\langle u,\tilde b\rangle+\int_{||y||\leq R}\big(e^{i\langle u,y\rangle}-1-i\langle u,y\rangle\ind_{||y||\leq 1}\big) \nu(dy)
  \end{eqnarray*}
  and $\tilde b=b-\int_{||y||>R} y\ind_{||y||\leq 1} \nu(dy)$. As in \cite{HM},   
\begin{align}  
\nonumber  
&\hskip -30pt 
P(f(X)-E[f(X)]\geq x)\\
\nonumber 
&\leq P(f(Y_R)-E[f(X)]\geq x)+P(Z_R\not =0)\\  
\nonumber  
&\leq P(f(Y_R)-E[f(Y_R)]\geq x-|E[f(X)]-E[f(Y_R)]|)+ P(Z_R\not =0)\\
\label{eq:split1}
&\leq P\left(f(Y_R)-E[f(Y_R)]\geq x-  E[||Z_R||]\right)+ P(Z_R\not =0)
\end{align}
where we use the fact that $f$ is a $1$-Lipschitz function. Also as in \cite{HM} (see \eqref{eq:ZR-Poisson}), we know that 
\begin{eqnarray}  
\nonumber  
P(Z_R\not=0)&\leq& 1-\exp\left(\nu(\|u\|\geq R)\right)\\  
\nonumber  
&=&\disp 1-\exp\left(-\int_{\|u\|\geq R} \nu(du)\right)\\  
\nonumber   
&=&\disp 1-\exp\left(-\int_{\s}\sigma(d\xi)\int_{\|r\xi\|\geq R}\frac{dr}{r^{1+\alpha}}\right)\\  
\label{eq:maj2}  
&=&\disp 1-\exp\left(-\frac{\s}{\alpha R^\alpha}\right)\\  
\label{eq:maj3}  
&\leq &\frac{\s}{\alpha R^\alpha}.
\end{eqnarray}  
In order to study the deviation of $f(Y_R)$, we shall apply in Sections \ref{sec:petit} and \ref{sec:intermediaire1} 
the lemma given below. This lemma generalizes, to arbitrary order $n\geq 2$, Lemma 2 in \cite{HM} which corresponds to $n=2$.

\begin{lemma}   
\label{lemme:bis3}  
Let $f:\rit^d\to \rit$ be a $1$-Lipschitz function and $Y_R\in \rit^d$ be a random vector whose L\'evy measure $\nu_R$ is the stable one truncated 
at $R>0$. Let $\delta>0$.

For any $n\geq 2$,  let $u_n(\alpha,\delta)$ be the unique non-zero (thus positive) solution of   
$$
\disp e^u-1-\frac{\delta(n-1)}{2-\alpha} u=0,   
$$   
and let   
\begin{equation}
\label{eq:un*}
u_n^*(\alpha,\delta)=\min_{1<k<\frac{n-1+\alpha}2}  \left(\left(k!\frac{\delta(k+1-\alpha)}{(n-k)(2-\alpha)}\right)^{1/(k-1)}\right)\wedge u_n(\alpha,\delta).  
\end{equation}
Then for all   
\begin{equation}  
\label{eq:Cond1}  
x\leq x_0:= (1+\delta(n-1))\frac{\s R^{1-\alpha}}{2-\alpha} u_n^*(\alpha,\delta),  
\end{equation}   
 we have   
\begin{equation}  
\label{eq:bis2-1}  
P(f(Y_R)-E[f(Y_R)] \geq x)\leq \exp\left(-\frac {(2-\alpha)x^2}{2(1+\delta(n-1))\s  R^{2-\alpha}}\right).  
\end{equation}  
\end{lemma}  
\begin{Rem}
{\rm
In the sequel, except stated otherwise, $\delta$ is mainly taken to be $1$ and $u_n(\alpha,\delta), u_n^*(\alpha,\delta)$ will simply be denoted 
$u_n(\alpha),u_n^*(\alpha)$. $\delta>0$ will be used in Proposition \ref{prop:D-s1.1} and Remark \ref{rem:petit2} recovering the Gaussian deviation from the stable one by letting $\delta\to 0$. 
}
\end{Rem}
\begin{Proof}  
 From  Theorem 1 in \cite{H02}, we know that   
\begin{equation}  
\label{eq:dev0}  
P(f(Y_R)-E[f(Y_R)]\geq x)\leq \exp\left(-\int_0^x h_R^{-1}(s) ds \right), \med x>0
\end{equation}  
with the function $h_R$ given by   
$$  
h_R(s) = \int_{\rit^d}  ||u||(e^{s ||u||}-1)\:\nu_R (dy), \quad s >0.
$$
Using for any $n \geq 2$,  
\begin{equation}  
\label{eq:h-bound2}  
e^{su}-1\leq \sum_{k=1}^{n-1} \frac{s^k}{k!} u^k + \frac{e^{sK}-\sum_{k=0}^{n-1}s^k K^k/k!}{K^n} u^n,   
 \qquad 0 \leq u \leq K, \ s\geq 0,   
\end{equation}   
 we have   
\begin{align}   
\nonumber  
& h_R(s)\leq  \int_{\rit^d}
\Big(\sum_{k=1}^{n-1} \frac{s^k}{k!}||u||^{k+1}\\  
\nonumber   
&\hskip 4cm +\frac{e^{s R}-\sum_{k=0}^{n-1}s^k R^k/k!}{R^n}||u||^{n+1}\Big) \nu_R(dy)   \\  
\nonumber  
&\hskip 1cm \leq \sum_{k=1}^{n-1} \frac{s^k}{k!}\int_{B(0,R)} ||u||^{k+1}\; \nu(dy)\\  
\nonumber  
&\hskip 4cm +\frac{e^{s R}-\sum_{k=0}^{n-1}s^k R^k/k!}{K^n} \int_{B(0,R)}  ||u||^{n+1} \; \nu(dy)\\
\nonumber  
&\hskip 1cm \leq\sum_{k=1}^{n-1}\left(\alpha_k-R^{k-n}\alpha_{n-1}\right)\frac{s^k}{k!}+\frac{\alpha_{n+1}}{R^n} (e^{s R}-1)\\
\nonumber
&\hskip 1cm \leq \sum_{k=1}^{n-1} \frac {(s R)^k}{k!}\frac{\s R^{1-\alpha} (n-k)}{(k+1-\alpha)(n+1-\alpha)}+\frac{\s R^{1-\alpha}}{n+1-\alpha} (e^{s R}-1)\\
\label{eq:alpha4}
&\hskip 1cm \leq \sum_{1\leq k<\left[\frac{n-1+\alpha}2\right]} \frac {(s R)^k}{k!}\frac{\s R^{1-\alpha} (n-k)}{(k+1-\alpha)(n+1-\alpha)}\\
\nonumber
&\hskip 5cm +\left(1+n-\left[\frac{n-1+\alpha}2\right]\right)\frac{\s R^{1-\alpha}}{n+1-\alpha} (e^{s R}-1),
\end{align}  
where we set $\alpha_k=\int_{B(0,R)}||u|| \; \nu(du)=\frac{\s}{k-\alpha}R^{k-\alpha}$ and where the last line follows from 
the fact that for    $k\geq \frac{n-1+\alpha}2$, the $k$-th summand in the above sum is bounded  by the last exponential term since $n-k\leq k+1-\alpha$.
  
For $1<k<\frac{n-1+\alpha}2$, set $u_k(\alpha)=\disp \left(k!\frac{\delta(n-1)(k+1-\alpha)}{(n-k)(2-\alpha)}\right)^{1/(k-1)}$   
and observe that for $s R\leq u_k(\alpha)$, the $k$-th term in the max is bounded by $\delta$ times the first one.   
Denote also by   
 $u_n(\alpha)$ the unique positive solution of   
$$  
e^u-1-\frac{\delta(n-1)}{2-\alpha} u=0.  
$$  
 Next, let  
$$  
u_n^*(\alpha)=   
 u_n(\alpha)   
 \wedge   
 \min_{1<k<\frac{n-1+\alpha}2} u_k(\alpha).  
$$  
For $s\leq \frac 1{ R} u_n^*(\alpha)$, all the terms in \eqref{eq:alpha4} are bounded by $\delta$ times the linear term, so that    
\begin{eqnarray*}  
h_R(s)&\leq& (1+\delta(n-1))\frac{\s R^{1-\alpha}}{n+1-\alpha} \frac{\delta(n-1)}{2-\alpha}s R\leq (1+\delta(n-1))\frac{\s R^{1-\alpha}}{2-\alpha} s R.   
\end{eqnarray*}  
 Hence, for $t\leq x_0:= (1+\delta(n-1))\frac{\s R^{1-\alpha}}{2-\alpha} u_n^*(\alpha)$, we can take   
$$  
h_R^{-1}(t)=\frac {(2-\alpha)t}{(1+\delta(n-1))\s R^{2-\alpha}},  
$$   
 which finally yields \eqref{eq:bis2-1} from \eqref{eq:dev0}. \CQFD  
\end{Proof}

  
 \begin{Rem}[on $u_k(\alpha)$, $u_n(\alpha)$, $u_n^*(\alpha)$, $x_0(n)$]
{\rm $ $

From now on and except stated otherwise, $\delta=1$.

\begin{itemize}   

\item The larger $n$, the worse the deviation \eqref{eq:bis2-1} becomes.  

\item For $n\leq 3$, the range   
 $1 < k < \frac{n-1+\alpha}2$, in the minimum defining $u_n^*(\alpha)$ is empty and thus $u_n^*(\alpha)=u_n(\alpha)$.  

\item As in (6.6) of \cite{BHP}, we have   
\begin{equation}  
\label{eq:u2}  
\log \frac{n-1}{2-\alpha}< u_n(\alpha)< 2\log \frac{n-1}{2-\alpha}.  
\end{equation}  
Moreover, for $n$ respectively larger than $5, 13, 18$, observe that we have  $u_n(\alpha)\geq 1,2,3$.  

\item   
 It is easy to see that for $k\geq 2$,   
$$  
k!\frac{n-1}{n-k}\frac{k+1-\alpha}{2-\alpha}\geq k!\frac{3-\alpha}{2-\alpha}\geq 3^{k-1},   
$$  
so that $u_k(\alpha)\geq 3$ and since moreover $\disp u_2(\alpha)=2\frac{n-1}{n-2}\frac{3-\alpha}{2-\alpha}\leq 4\frac{3-\alpha}{2-\alpha}$, for $n\geq 3$, we have   
$$  
n_0\leq u_n^*(\alpha)\leq 4\frac{3-\alpha}{2-\alpha},  
$$  
and  we can take  
$$  
x_0 = n n_0 \frac{\s  R^{1-\alpha}}{2-\alpha},   
$$   
with $n_0=1,2,3$, for respectively $n\geq 5,13,18$.  

\item When $n$ is fixed and $\alpha \nearrow 2$, we have for two constants $K_1(n)$ and $K_2(n)$  
$$  
\frac{K_1(n)}{(2-\alpha)^{2/(n+1)}}<\min_{1<k<\frac{n-1+\alpha}2} u_k(\alpha) \leq \frac{K_2(n)}{2-\alpha},   
$$   
so that, from \eqref{eq:u2}, $u_n^*(\alpha)=u_n(\alpha)$ is at most of order $\log\frac{n-1}{2-\alpha}$.  

\end{itemize}  
}
\end{Rem}


\subsection{Lower range for the stable deviation}
\label{sec:petit}
The deviation \eqref{eq:dev-def} for small values of $x$ is given by the following:
\begin{prop}  
\label{prop:dev-petit2}  
Let $X$ be a stable random vector in $\rit^d$ of index $\alpha\in (1,2)$ and let $f:\rit^d\to\rit$ be a $1$-Lipschitz function.   
Then for all $n$ large enough and for all
\begin{equation}
\label{eq:lambda11}
\lambda \in \left]\frac{\s}{\alpha-1}\left(1+\sqrt{\frac{2(\alpha-1)^2n \s^{\frac{2-\alpha}{\alpha-1}}}{\alpha(2-\alpha)}}\right), \frac \s{\alpha-1}\left(1+\frac{\alpha-1}{2-\alpha}n u_n^*(\alpha)\right)\right[,
\end{equation}
there exists $x_1(n,\alpha, \lambda)>0$ such that for all $0\leq x\leq  x_1(n,\alpha,\lambda)$, 
\begin{align}
\nonumber
&P(f(X)-E[f(X)]\geq x)\\
\label{eq:petit1}
&\leq 
 \exp\left(-\frac{2-\alpha}{2n\s^{1/(\alpha-1)}}(\lambda-\frac \s{\alpha-1})^2
\left(\frac x\lambda \right)^{\frac \alpha{\alpha-1}} \right)
+\frac \s\alpha \left(\frac x\lambda \right)^{\frac \alpha{\alpha-1}}\leq 1.
\end{align}
\end{prop}  

In \eqref{eq:petit1}, our purpose is to investigate the order of deviation from the mean of Lipschitz function of stable random vector 
and for $x$ small. The order obtained is essentially $\exp\left(-cx^{\alpha/(\alpha-1)}\right)$ for some explicit constant $c$. 
Note that the exponent $\frac \alpha{\alpha-1}$ of $x$ in this bound goes to $2$ when $\alpha$ goes to $2$, 
this is reminiscent of the Gaussian case. A more precise statement is given in Remark \ref{rem:petit2} in connexion with Gaussian bound of deviation.  

\begin{Proof}
In order to investigate deviation for small values of $x$, in this section the level of truncation is chosen by setting 
$$
x= \frac{\lambda}{R^{\alpha-1}}
$$
where $\lambda>0$ is as in the Proposition \ref{prop:dev-petit2}. 

Set also $u(R)= (\lambda-\frac{\sigma}{\alpha-1})/R^{\alpha-1}=\tilde \lambda /R^{\alpha-1}$. Note that $Z_R$ has a compound Poisson structure and is the same in law as $Z_R=\sum_{k=0}^{N}Z_k$,
where $Z_0=0$, and the $Z_k$, $k\geq 1$, are i.i.d. random vectors with law $\disp \frac{\nu_{Z_R}}{\nu\{||u||>R\}}$ and $N$ is an 
independent Poisson random variable with intensity $\nu\{||u||>R\}$. Hence, for any $B\in{\cal B}(\rit^d)$,
\begin{equation}
\label{eq:ZR-Poisson}
P(Z_R\in B)=\sum_{n=0}^{+\infty}\frac{e^{-\nu\{||u||>R\}}}{n!}(\nu\{||u||>R\})^n P\left(\sum_{k=0}^{n}Z_k\in B\right).
\end{equation}
Thus,  
\begin{eqnarray}
\nonumber
E[||Z_R||]&\leq&
\sum_{n=0}^{+\infty}\frac{e^{-\nu\{||u||>R\}}}{n!}(\nu\{||u||>R\})^n E\left[\sum_{k=0}^{n}||Z_k||\right]\\
\nonumber
&=&E[||Z_1||] e^{-\nu\{||u||>R\}}\sum_{n=1}^{+\infty}\frac{\nu\{||u||>R\}^n}{(n-1)!}\\
\nonumber 
&=&E[||Z_1||] \: \nu\{||u||>R\}\\
\label{eq:ZR1}
&=&\int_{||u||>R}||u||\nu(du)=\frac{\s}{\alpha-1} R^{1-\alpha}.
\end{eqnarray}
Moreover, we also have a lower bound as follows:
\begin{eqnarray}
\nonumber
E[||Z_R||]&=&\int_0^{+\infty} P(||Z_R||>x) dx\\
\nonumber
&=&\sum_{n=0}^{+\infty}\frac{e^{-\nu\{||u||>R\}}}{n!}(\nu\{||u||>R\})^n \int_0^{+\infty}P\left(\left\|\sum_{k=0}^{n}Z_k\right\|>x\right)dx\\
\nonumber
&\geq&e^{-\nu\{||u||>R\}}\nu\{||u||>R\}\int_0^{+\infty}P(||Z_1||>x)dx\\
\nonumber 
&=& e^{-\nu\{||u||>R\}}\nu\{||u||>R\}E[||Z_1||]\\
\label{eq:ZR2}
&=&\frac{\s}{(\alpha-1)R^{\alpha-1}}e^{-\frac{\s}{\alpha R^\alpha}}.
\end{eqnarray}
From \eqref{eq:ZR1} and \eqref{eq:ZR2}, we get
\begin{equation}
\label{eq:EZR}
\frac{\s}{(\alpha-1)R^{\alpha-1}}e^{-\frac{\s}{\alpha R^\alpha}}
\leq E[||Z_R||]\leq 
\frac{\s}{(\alpha-1) R^{\alpha-1}}.
\end{equation}
Next, we go back to the study of
\begin{equation}
\label{eq:2}  
P(f(X)-E[f(X)]\geq x)\leq P(f(Y_R)-E[f(Y_R)]\geq u(R))+ P(Z_R\not =0).
\end{equation}

\noindent
For the second summand, a bound is given in \eqref{eq:maj3}. 
For the first summand, applying Lemma \ref{lemme:bis3}, we have 
\begin{equation}
\label{eq:dev1R}
P(f(Y_R)-E[f(Y_R)]\geq u(R))\leq \exp\left(-\frac{2-\alpha}{2n\s}\frac{u(R)^2}{R^{2-\alpha}}\right)
\end{equation}
as long as $0<u(R)<n\frac{\s R^{1-\alpha}}{2-\alpha}u_n^*(\alpha)$ that is as long as
$$
\left(\lambda-\frac \s{\alpha-1}\right)\frac 1{R^{\alpha-1}}
<\frac{n\s u_n^*(\alpha)}{2-\alpha}R^{1-\alpha}
$$
that is for 
\begin{equation}
\label{eq:lambda1}
\lambda<\frac \s{\alpha-1}\left(1+\frac{\alpha-1}{2-\alpha}n u_n^*(\alpha)\right). 
\end{equation} 
From \eqref{eq:2} and \eqref{eq:dev1R}, we have 
 \begin{align}
\nonumber
&P(f(X)-E[f(X)]\geq x)\\
\label{eq:Cor}
&\leq 
 \exp\left(-\frac{2-\alpha}{2n\s^{1/(\alpha-1)}}(\lambda-\frac \s{\alpha-1})^2
\left(\frac x\lambda \right)^{\frac \alpha{\alpha-1}} \right)
+\frac \s\alpha \left(\frac x\lambda \right)^{\frac \alpha{\alpha-1}}.
\end{align}
The bound \eqref{eq:Cor} makes sense if the right-hand side is smaller than $1$, this is true at $0^+$ if 
\begin{equation}
\label{eq:lambda2}
\lambda> \frac{\s}{\alpha-1}+\sqrt{\frac{2n \s^{\frac{\alpha}{\alpha-1}}}{\alpha(2-\alpha)}},
\end{equation} 
and \eqref{eq:lambda1}, \eqref{eq:lambda2} are compatible if 
\begin{equation}
\label{eq:rem}
\frac{2(2-\alpha)\s^{\frac{2-\alpha}{\alpha-1}}}{\alpha n}<u_n^*(\alpha)^2, 
\end{equation}
and this is true if $n$ is large enough or if $\alpha$ is close enough to $2$. \CQFD
\end{Proof}

\begin{Rem}[Comments on the bound and on the range of Proposition \ref{prop:dev-petit2}]
\label{rem:petit1}
{\rm $ $

$\bullet$ In fact, the bound \eqref{eq:petit1} holds for all $x>0$ but makes sense only for $x<x_1(n,\alpha, \lambda)$. 

$\bullet$ The value $x_1(n,\alpha,\lambda)$ in the bound of the domain of Proposition \ref{prop:dev-petit2} can be made more precise: $\disp x_1(n,\alpha,\lambda)=\lambda u_0(n,\alpha,\lambda)^{\frac{\alpha-1}{\alpha}}$ where $u_0(n,\alpha,\lambda)$ is the unique positive solution of 
$$
 \exp\left(-\frac{2-\alpha}{2n\s^{1/(\alpha-1)}}(\lambda-\frac \s{\alpha-1})^2
u \right)
+\frac \s\alpha u=1.
$$
In particular, we have 
$$
u_0(n,\alpha,\lambda)\geq \frac{2n\s^{1/(\alpha-1)}}{(2-\alpha)(\lambda-\frac \s{\alpha-1})^2}\ln\left(\frac{\alpha(2-\alpha)}{2n\s^{\alpha/(\alpha-1)}}\big(\lambda-\frac \s{\alpha-1}\big)^2\right)
$$
so that we can take 
\begin{equation}
\label{eq:x1}
x_1(n,\alpha,\lambda)=\lambda \s^{1/\alpha}\left(\frac{2n}{(2-\alpha)(\lambda-\frac \s{\alpha-1})^2}\ln\left(\frac{\alpha(2-\alpha)}{2n\s^{\alpha/(\alpha-1)}}\big(\lambda-\frac \s{\alpha-1}\big)^2\right)\right)^{\frac{\alpha-1}\alpha}.
\end{equation}

$\bullet$ In order to optimize the bound \eqref{eq:petit1} with repect to $\lambda$, observe that the first summand (which gives the right order) gives the best bound for  
$$
\lambda=\lambda_0(\alpha):=\frac{\alpha\s}{(2-\alpha)(\alpha-1)}. 
$$
Denote by $]\lambda_1(n,\alpha),\lambda_1(n,\alpha) [$, the interval in \eqref{eq:lambda11}.\\
We have  $\lambda_0(\alpha)\geq \lambda_1(n,\alpha)$ if $n\s^{\frac{2-\alpha}{\alpha-1}}\leq \frac{2\alpha}{\alpha-1}$ and $\lambda_0(\alpha)\leq \lambda_2(n,\alpha)$ if $nu_n^*(\alpha)\geq 2$ so that usually 
$$
\lambda_0(\alpha)\leq \lambda_1(n,\alpha)
$$
and the best order should thus occurs for  $\lambda=\lambda_1(n,\alpha)$. But for this choice $x_1(n,\alpha,\lambda)=0$ and the domain of deviation is empty. There is thus no optimal choice for $\lambda\in ]\lambda_1(n,\alpha),\lambda_1(n,\alpha) [$: for $\lambda=\lambda_1(n,\alpha)$, the deviation is the best but the domain is empty, while for $\lambda=\lambda_2(n,\alpha)$ the domain is the largest but the deviation is the worst. See below in Remark \ref{rem:petit3} for a deviation without extra parameter. 

$\bullet$ In \eqref{eq:petit1}, the exponential is the main term of the bound. For any $\varepsilon>0$, we thus can rewrite \eqref{eq:petit1},  for $x$ small enough
$$
P(f(X)-E[f(X)]\geq x)\leq (1+\varepsilon)
 \exp\left(-\frac{2-\alpha}{2n\s^{1/(\alpha-1)}}(\lambda-\frac \s{\alpha-1})^2
\left(\frac x\lambda \right)^{\frac \alpha{\alpha-1}} \right). 
$$
In fact, it gives the order of stable deviation for small values of $x$: roughly speaking, it is in $\disp \exp\big(-cx^{\frac{\alpha}{\alpha-1}}\big)$. 
}
\end{Rem}

\begin{Rem}[Gaussian deviation]
\label{rem:petit2}
{\rm 
We  can recover the Gaussian deviation from the bound \eqref{eq:petit1} letting $\alpha$ goes to $2$ properly.
To this way and following \cite{M2}, consider a stable random vector $X^{(\alpha)}$ whose L\'evy measure has spherical component $\sigma$ given by the sum of Dirac measures at the points
$(0,\dots, 0,\pm 1,0, \dots,0)$ and with total mass $\s =2-\alpha$. The components of $X^{(\alpha)}$ are thus independent and when $\alpha$ goes to $2$, the vector $X^{(\alpha)}$ converges in distribution to a standard Gaussian random vector~$W$. Moreover, note that when we take limit $\alpha \to 2$ with $\s=2-\alpha$ in \eqref{eq:lambda1} and \eqref{eq:x1}, the ranges for $\lambda$  and for $x$ in Proposition \ref{prop:dev-petit2} becomes $\lambda\in(0,\infty)$ and $x\in (0, \infty)$ while the bound \eqref{eq:petit1} becomes $\exp(- x^2/(2n))$. This is not exactly the classical bound for  Gaussian deviation. But, first note that we could replace $n$ large in Proprosition \ref{prop:dev-petit2} by $\alpha$ close enough to $2$ (see inequality \eqref{eq:rem} in the proof). Then  apply Lemma \ref{lemme:bis3} with arbirary $\delta>0$ and take limit $\alpha\to 2$ in \eqref{eq:petit1}, we obtain similarly 
$$
P(f(W)-E[f(W)]\geq x)\leq \exp\left(-\frac {x^2}{2n_\delta}\right),
$$
for any $x>0$. Finally, letting $\delta\to 0$ yields the Gaussian deviation bound for all $x>0$:
$$
P(f(W)-E[f(W)]\geq x)\leq \exp\left(-\frac {x^2}{2}\right).
$$

The same can be made from the intermediate regime of deviation studied in Section \ref{sec:intermediaire1}, see the Remark \ref{rem:Gauss}.

}\end{Rem}

\begin{Rem}
\label{rem:petit3}
{\rm
We can give a deviation bound for small values of $x$ without introducing an extra parameter $\lambda$ as in Proposition \ref{prop:dev-petit2}. To this way, let $\varepsilon>0$. For all $n\geq 5$ (or for $\alpha$ close enough to $2$), there exists $x_0(n,\varepsilon)>0$ such that for all $0\leq x< x_0(n,\varepsilon)$, 
\begin{equation}
\label{eq:dev-petit}
P(f(X)-E[f(X)]\geq x)\leq (1+\varepsilon) \exp\left(-\frac{2-\alpha}{2n\s^{1/(\alpha-1)}}
\left(\frac{\alpha-1}{\alpha}\right)^{\frac \alpha{\alpha-1}} x^{\frac \alpha{\alpha-1}}\right).
\end{equation}

The drawback of this bound is its range since we do not know explicitly $x_0(n, \varepsilon)$. 

The proof follows the same lines of reasoning as that of Proposition \ref{prop:dev-petit2} but with the level of truncation chosen by setting 
\begin{equation*}
x= \alpha E[||Z_R||].
\end{equation*}
Set also $u(R)= (\alpha-1) E[||Z_R||]$ and study the summands of \eqref{eq:2}. Applying Lemma \ref{lemme:bis3} and using 
\begin{eqnarray*}
\frac{u(R)^2}{R^{2-\alpha}}&=&\left(\frac{\alpha-1}{\alpha}\right)^2\frac{x^2}{R^{2-\alpha}}
\geq\left(\frac 1{\s}\right)^{\frac{2-\alpha}{\alpha-1}} \left(\frac{\alpha-1}{ \alpha}\right)^{\frac \alpha{\alpha-1}}x^{\frac \alpha{\alpha-1}},
\end{eqnarray*}
the first summand in the right hand side of \eqref{eq:2} is bounded as follows
\begin{equation*}
\label{eq:dev2R}
P(f(Y_R)-E[f(Y_R)]\geq u(R))\leq \exp\left(-\frac{2-\alpha}{2n\s^{1/(\alpha-1)}}
\left(\frac{\alpha-1}{\alpha}\right)^{\frac \alpha{\alpha-1}} x^{\frac \alpha{\alpha-1}}\right)
\end{equation*}
as long as $0<\frac{\alpha-1}\alpha x=u(R)<n\frac{\s R^{1-\alpha}}{2-\alpha}u_n^*(\alpha)$ that is as long as
$$
E[||Z_R||]<\frac{n\s u_n^*(\alpha)}{(2-\alpha)(\alpha-1)}R^{1-\alpha},
$$
which, from the right-hand side of \eqref{eq:EZR}, is true if $\disp\frac{\s}{\alpha-1}<\frac{n\s u_n^*(\alpha)}{(2-\alpha)(\alpha-1)}$, 
that is if $\disp n u_n^*(\alpha)>2-\alpha$, which is true at least for $\alpha$ close enough to $2$ or for $n\geq 5$ since then $u_n^*(\alpha)\geq 1$.
Using \eqref{eq:maj3} for the second summand in the right-hand side of \eqref{eq:2}, we derive 
\begin{equation}
\label{eq:dev-petit3} 
P(f(X)-E[f(X)]\geq x)
\leq \exp\left(-\frac{2-\alpha}{2n\s^{1/(\alpha-1)}}\left(\frac{\alpha-1}{\alpha}\right)^{\frac \alpha{\alpha-1}} x^{\frac \alpha{\alpha-1}}\right)+\frac{\s}{\alpha R^\alpha}. 
\end{equation}
But from \eqref{eq:EZR},  when $x$ goes to $0$, we have $R\to+\infty$ or $R\to 0$. 
Moreover, when $x\to 0$, $E[||Z_R||]\to 0$, so that $Z_R\to 0$ in $L^1(\rit^d)$ and $Z_R\Rightarrow \delta_0$. 
This implies the convergence of $\nu_{Z_R}=\nu_{B(0,R)^c}$ to $0$, from which we derive $R\to +\infty$. 
Thus when $x$ is small, the main term in the bound \eqref{eq:dev-petit3} in given by the exponential. Since the second term goes to $0$, for any $\varepsilon>0$, there is some $x_0(n,\varepsilon)$ such that for $x\leq x_0(n,\varepsilon)$, the bound \eqref{eq:dev-petit} holds.
}
\end{Rem}


\subsection{Intermediate range for the stable deviation}
\label{sec:intermediaire1}
  
In this section, we study the deviation \eqref{eq:dev-def} for intermediate values of $x$. To this end, the level of truncation is changed in \eqref{eq:split1}. 
Moreover in order to derive the Gaussian deviation from the stable one as a limiting case, we shall use Lemma 
\ref{lemme:bis3} with parameter $\delta>0$. In the limiting case, $\delta$ will goes to $0$. First, we can state:
\begin{prop}  
\label{prop:D-s1.1}  
Let $X$ be a stable random vector in $\rit^d$ of index $\alpha\in (1,2)$ and let $f:\rit^d\to\rit$ be a $1$-Lipschitz function.   
Let $\delta>0$ and $n\in\nit$, $n\geq 2$. Set $n_\delta=(1+\delta(n-1))$. Then for any $\varepsilon$ satisfying 
\begin{equation}  
\label{eq:varep1}  
\varepsilon> \frac{(2-\alpha)e}{2\alpha n_\delta},  
\end{equation}  
we have
\begin{equation}  
\label{eq:Deviation1}  
P(f(X)-E[f(X)]\geq x)\leq (1+\varepsilon)\exp\left(-\frac {2-\alpha}{2n_\delta\s } 
\left(\frac{x}{1+\frac {2-\alpha}{2n_\delta(\alpha-1)u_1(n_\delta,\alpha)}}\right)^{\alpha}\right),
\end{equation}   
for any $x$ such that   
\begin{align}  
\label{eq:range1-x}  
&\frac {2n_\delta \s }{2-\alpha}u_1(n_\delta,\alpha)\left(1+\frac{2-\alpha}{2n_\delta(\alpha-1) u_1(n_\delta,\alpha)}\right)^\alpha<x^\alpha<\\  
\nonumber   
&  
\frac {2n_\delta \s }{2-\alpha}\big(u_2(n_\delta,\alpha)\wedge u_n^*(\alpha,\delta)/2\big)\left(1+\frac{2-\alpha}{2n (\alpha-1)\big(u_2(n_\delta,\alpha)\wedge u_n^*(\alpha,\delta)/2\big)}\right)^\alpha
\end{align}  
where $u_1(n_\delta,\alpha)$ and $u_2(n_\delta,\alpha)$ are the solutions of 
$$
e^{u}-\frac{2 n_\delta\alpha\varepsilon}{2-\alpha}u=0
$$
and $u_n^*(\alpha,\delta)$ is given by \eqref{eq:un*}.
\end{prop}

\begin{Proof}
First, using the right-hand side of \eqref{eq:EZR}, the bound \eqref{eq:split1} becomes 
$$
P(f(X)-E[f(X)]\geq x)\leq P\left(f(Y_R)-E[f(Y_R)]\geq x-\frac{\s}{\alpha-1} R^{1-\alpha}\right)+ P(Z_R\not =0).
$$
Next choose the level of truncation by setting 
\begin{equation}  
\label{eq:Rx1}  
 R\left(1+\frac{\s }{(\alpha-1) R^\alpha}\right)=x.   
\end{equation}  
We thus have
\begin{equation}
\label{eq:1}  
P(f(X)-E[f(X)]\geq x)\leq P(f(Y_R)-E[f(Y_R)]\geq  R)+ P(Z_R\not =0).
\end{equation}

Applying Lemma \ref{lemme:bis3} with $x= R$ in \eqref{eq:bis2-1} to estimate the first summand in the right-hand side of \eqref{eq:1} 
and using \eqref{eq:maj3} for the second one, the bound \eqref{eq:1} then becomes    
\begin{equation}  
\label{eq:D}  
P(f(X)-E[f(X)]\geq x)\leq   
\exp\left(-\frac{(2-\alpha)R^\alpha}{2n_\delta\s }\right)+\frac{\s}{\alpha R^\alpha}.  
\end{equation}  
as long as, using \eqref{eq:Cond1} with  $x= R$,  
\begin{equation}  
\label{eq:R-1}  
R^\alpha \leq n_\delta\frac{\s}{2-\alpha} u_n^*(\alpha,\delta).  
\end{equation}  

Next, we compare the two summands in the right-hand side of \eqref{eq:D}.  
To do this, set $u=\frac{(2-\alpha)R^\alpha}{2n_\delta\s }$ and let us compare $\frac{2-\alpha}{2n \alpha u}$ to $e^{-u}$ 
studying the function   
\begin{equation}  
\label{eq:h4}  
h_{n,\alpha}(u)=e^{u}-\frac{2n_\delta \alpha\varepsilon}{2-\alpha}u,   
\end{equation}  
Note that $h_{n,\alpha}$ has a unique minimum at $u_0=\log \frac {2n_\delta \alpha \varepsilon}{2-\alpha}$, which is negative because of \eqref{eq:varep1}. 
 We thus have   
\begin{equation}  
\label{eq:21}  
P(f(X)-E[f(X)]\geq x)\leq (1+\varepsilon)\exp\left(-\frac{2-\alpha}{2n_\delta \s }R^\alpha\right),  
\end{equation}  
for   
\begin{equation}  
\label{eq:domaine-R}  
\frac {2n_\delta \s}{2-\alpha}u_1(n_\delta,\alpha)< R^\alpha<  \frac {2n_\delta \s}{2-\alpha}u_2(n_\delta,\alpha),  
\end{equation}  
and still under the condition \eqref{eq:R-1}, with moreover $R$ given by the relation \eqref{eq:Rx1}.   
  
\vskip 10pt  
  
We now express \eqref{eq:21} and its conditions \eqref{eq:R-1} and \eqref{eq:domaine-R}  in terms of $x$. Since from \eqref{eq:Rx1}  
\begin{equation}  
\label{eq:R-x1}  
 R\leq x\leq  R\left(1+\frac {2-\alpha}{2n_\delta(\alpha-1) u_1(n_\delta,\alpha)}\right),  
\end{equation}  
the deviation bound \eqref{eq:21} becomes \eqref{eq:Deviation1}.

To express the range in terms of $x$, note that we have 
$\disp x =\theta(R^\alpha)$ with   
\begin{equation*}  
\label{eq:theta}  
\theta(u)=u^{1/\alpha}\left(1+\frac \s{(\alpha-1)u}\right).  
\end{equation*}  
The function $\theta$ is an increasing bijection from $[\s, +\infty)$ to $[\frac  \alpha{\alpha-1} \s^{1/\alpha}, +\infty)$.   
  
Then, note that $\s\leq \frac{2n_\delta\s}{2-\alpha}u_1(n_\delta,\alpha)$, that is $\frac{2-\alpha}{2n_\delta}\leq u_1(n_\delta,\alpha)$, which is equivalent to have $h_{n,\alpha}(\frac{2-\alpha}{2n_\delta})\geq 0$ and $h_{n,\alpha}'(\frac{2-\alpha}{2n_\delta})\leq 0$. 
Indeed, writing $\varepsilon =\frac{(2-\alpha)\tilde e}{2n_\delta\alpha}$ in \eqref{eq:varep1}, with some $\tilde e>e$, we have    
\begin{eqnarray*}  
h_{n,\alpha}\left(\frac{2-\alpha}{2n_\delta}\right)&=&\exp\left(\frac{2-\alpha}{2n_\delta}\right)-\alpha \varepsilon=\exp\left(\frac{2-\alpha}{2n_\delta}\right)-\frac{2-\alpha}{2n_\delta}\tilde e.  
\end{eqnarray*}  
Let $\eta(u)=e^u-\tilde e u$. Since $\frac{2-\alpha}{2n_\delta} \leq 1/2<1$, we have $\eta(\frac{2-\alpha}{2n_\delta})>0$ as long as $\tilde e$ is such that $\tilde e\leq 2\sqrt e$. 
We dedude $h_{n,\alpha}(\frac{2-\alpha}{2n_\delta})\geq 0$.  
  
Next   
$\disp  
h'_{n,\alpha}\left(\frac{2-\alpha}{2n_\delta}\right)=\exp\left(\frac{2-\alpha}{2n_\delta}\right)-\frac{2n_\delta\alpha}{2-\alpha}\varepsilon  
\leq \exp\left(\frac{2-\alpha}{2n_\delta}\right)-e  
$,  
which is non-positive since $2-\alpha\leq 2n_\delta$. Finally,   
$$  
\s\leq \frac{2n_\delta\s}{2-\alpha}u_1(n_\delta,\alpha).   
$$  
  
Using $\theta$ in the domain where it defines an increasing bijection, the conditions  \eqref{eq:R-1} and \eqref{eq:domaine-R} can be rewritten as:  
$$  
\theta\left(\frac {2n_\delta \s }{2-\alpha}u_1(n_\delta,\alpha)\right)  
< x<  
 \theta\left(\frac {2n_\delta \s }{2-\alpha}\big(u_2(n_\delta,\alpha)\wedge u_n^*(\alpha,\delta)/2\big)\right).  
$$  
The range of the deviation \eqref{eq:Deviation1} is thus   
\eqref{eq:range1-x} when $\varepsilon$ satisfies  
$$  
\frac{(2-\alpha)e}{2n_\delta\alpha}<\varepsilon<\frac{(2-\alpha)\sqrt e}{n_\delta\alpha }.  
$$  
Finally, it is clear than if the deviation \eqref{eq:Deviation1} holds for some $\varepsilon>0$, it also holds for any larger $\varepsilon$. 
The condition on $\varepsilon$ can thus be reduced to \eqref{eq:varep1}. \CQFD
\end{Proof}  
  
  \vskip 5pt  
  
\begin{Rem} 
{\rm  $ $

$\bullet$  From \eqref{eq:varep1}, note that   
$$  
h_{n,\alpha}\left(\log\left( \frac {2n_\delta \alpha \varepsilon}{2-\alpha}\right)\right)=\frac {2n_\delta \alpha \varepsilon}{2-\alpha}\left(1-\log\left( \frac {2n_\delta \alpha \varepsilon}{2-\alpha}\right)\right)<0  
$$  
so that $\disp u_2(n_\delta,\alpha)> \log\left( \frac {2n_\delta\alpha\varepsilon}{2-\alpha}\right)>1$.  

$\bullet$ The parameter $n\geq 2$ is free in Proposition \ref{prop:D-s1.1}. Thus hanging $n$ allows to shift the range \eqref{eq:range1-x} so that we can derive a $e^{-cx^\alpha}$ -type bound of deviation on a larger domain, by eventually changing $n$ and the generic constant $c$.  

}   
\end{Rem}

  
\begin{Rem}[Comments on the bounds in the domain \eqref{eq:range1-x}]  
{\rm   $ $
  
It is interesting to investigate the behavior of $\disp \frac{2n_\delta\alpha}{2-\alpha} u_i(n_\delta,\alpha)$, $i=1,2$, 
when $\alpha$ goes to $2$ or when $n$ goes to $+\infty$ in the domain given by \eqref{eq:range1-x}.   
  
To simplify notation, \eqref{eq:h4} is written as $h_a(u)=e^u-\frac{\varepsilon}{a}u$ with $\varepsilon>ae$ , and without loss of generality, 
we study the behavior of the zeros $u_1(a)<u_2(a)$ of $h_a(u)$ as $a\to 0$.   
  
Since $h_a(\log(\varepsilon/a))<0$ and $\lim_{a\to 0}h_a(a^\lambda)=1-\lambda\varepsilon$, for $\lambda<\varepsilon<\tilde \lambda$, we have   
$$  
\lambda a\leq  u_1(a)\leq \tilde \lambda a \leq \log(\varepsilon/a)\leq u_2(a).  
$$  
Thus,   
$$  
\lim_{a\to 0} \frac{u_1(a)}a=\varepsilon, \med \lim_{a\to 0} \frac{u_2(a)}a=+\infty,  
$$   
that is in our setting
\begin{equation}  
\label{eq:behavior}  
\lim_{\begin{subarray}{c}\alpha\to 2 \mbox{ {\small or} } n\to+\infty\end{subarray}}  
\frac{2n_\delta\alpha}{2-\alpha} u_1(n_\delta,\alpha)=1/\varepsilon  
,\med   
\lim_{\begin{subarray}{c}\alpha\to 2 \mbox{ {\small or} } n\to+\infty\end{subarray}}  
\frac{2n_\delta\alpha}{2-\alpha} u_2(n_\delta,\alpha)=+\infty.  
\end{equation}  
For $\alpha$ close to $2$ or for $n$ large enough, $\disp 1+\frac {2-\alpha}{2n_\delta(\alpha-1)u_1(n_\delta,\alpha)}\to 1+\frac{\alpha\varepsilon}{\alpha-1}$,  
so that for $\varepsilon<1$ and for $\alpha$ close to $2$ or for $n$ large enough, $\disp 1+\frac {2-\alpha}{2n_\delta(\alpha-1)u_1(n_\delta,\alpha)}\leq 1+\frac{\alpha}{\alpha-1}$. In this case, the deviation \eqref{eq:Deviation1} becomes   
\begin{equation*}  
\label{eq:Deviation1bis}  
P(f(X)-E[f(X)]\geq x)\leq (1+\varepsilon)\exp\left(-\frac {2-\alpha}{2n_\delta\s} \left(\frac{x}{1+\frac{\alpha}{\alpha-1}}\right)^{\alpha}\right),  
\end{equation*}   
and the range \eqref{eq:range1-x} reduces to 
\begin{equation*}  
\left(1+\frac {\alpha\varepsilon}{\alpha-1}\right)^{\alpha}\frac {\s}{\alpha\varepsilon}< x^\alpha<  
\left\{  
\begin{array}{ll}  
\disp \frac {2n_\delta\s}{2-\alpha}\log\left(\sqrt\frac{n_\delta-1}{2-\alpha}\right)&\mbox{ if } \alpha\to 2,\\  
\disp \frac {3n_\delta\s}{2-\alpha} &\mbox{ if } n\to+\infty.  
\end{array}  
\right .  
\end{equation*}  
}
\end{Rem}
\vskip 5pt

\begin{Rem}[More on Gaussian deviation]
\label{rem:Gauss}
$ $
{\rm 

$\bullet$ When $\alpha$ is close to $2$, we recover (i) of Theorem 2 in \cite{M2} and in particular still following \cite{M2}, we recover Gaussian deviation. To this way, as in Remark \ref{rem:petit2}, consider a stable random vector $X^{(\alpha)}$ whose L\'evy measure has spherical component $\sigma$ given by the sum of Dirac measures at the points  $(0,\dots, 0,\pm 1,0, \dots,0)$ and with total mass $\s =2-\alpha$. When $\alpha\to 2$, the vector $X^{(\alpha)}$ converges in distribution to a standard Gaussian random vector~$W$, so that \eqref{eq:Deviation1} yields 
$$
P(f(W)-E[f(W)]\geq x)\leq (1+\varepsilon)\exp\left(-\frac 1{2n_\delta} 
\left(\frac{x}{1+\varepsilon}\right)^2\right),
$$
for any $x>0$, since the left-hand side of \eqref{eq:range1-x} goes to $0$ when $\alpha$ goes to $2$ ($u_1(n_\delta,\alpha)\to 0$ as $\alpha\to 2$)
while its right-hand side goes to $+\infty$ ($u_2(n_\delta,\alpha)\to +\infty$ and $u_n^*(\alpha,\delta)\to+\infty$ as $\alpha\to 2$).

Finally, letting $\varepsilon\to 0$ and $\delta\to 0$ yields the Gaussian deviation bound for all $x>0$:
$$
P(f(W)-E[f(W)]\geq x)\leq \exp\left(-\frac {x^2}{2}\right).
$$

$\bullet$ Note also that letting $\delta\to 0$ in the range \eqref{eq:range1-x} does not allow to recover a deviation result for arbitrary small values of $x$. Proposition \ref{prop:dev-petit2} cannot thus be derived from Proposition~\ref{prop:D-s1.1} and vice versa. 
}
\end{Rem}
  

\subsection{Another intermediate range for stable deviation}  
\label{sec:intermediaire2}  
  
The deviation and the range obtained in Section \ref{sec:intermediaire1} depend on $u_1(n_\delta,\alpha)$ which is not explicit. 
In this section, using directly \eqref{eq:maj2} rather than the bound \eqref{eq:maj3}, we obtain the same type of result, 
probably less sharp,  but with more explicit bounds. 

\begin{prop}  
\label{prop:D-s2.1}  
Let $X$ be a stable random vector in $\rit^d$ of index $\alpha\in (1,2)$ and let $f:\rit^d\to\rit$ be a $1$-Lipschitz function.   
Then for any $\varepsilon$ satisfying
\begin{equation}  
\label{eq:varep11}  
\varepsilon\geq ec =\frac{(2-\alpha)e}{2n\alpha},  
\end{equation}  
we have
\begin{equation}  
\label{eq:Deviation2}  
P(f(X)-E[f(X)]\geq x)\leq (1+\varepsilon)\exp\left(-\frac {2-\alpha}{2n\s} \left(\frac{x}{1+\frac {2-\alpha}{2n(\alpha-1)\left(1-\left(\frac {2-\alpha}{2n\alpha}\right)^{2/3}\right)}}\right)^{\alpha}\right).  
\end{equation} 
for all $x$ such that 
\begin{align}  
\nonumber   
&\frac{2n\s}{2-\alpha}\left(1-\left(\frac {2-\alpha}{2n\alpha}\right)^{2/3}\right)\left(1+\frac{2-\alpha}{2n(\alpha-1)\left(1-\left(\frac {2-\alpha}{2n\alpha}\right)^{2/3}\right)}\right)^\alpha \\  
\label{eq:range-x2}  
&  
<x^\alpha<\frac{2n\s}{2-\alpha}\left(\left(1+\left(\left(\frac {2-\alpha}{2n\alpha}\right)^{2/3}\wedge 0.68\right)\right)\wedge u_n^*(\alpha)/2\right)\times\\  
\nonumber   
&\hskip 10pt   
\left(1+\frac{2-\alpha}{2n(\alpha-1)\left(\left(1+\left(\left(\frac {2-\alpha}{2n\alpha}\right)^{2/3}\wedge 0.68\right)\right)\wedge u_n^*(\alpha)/2\right)}\right)^\alpha  
\end{align}  
where $u_n^*(\alpha)$ is given in \eqref{eq:un*} with $\delta=1$.
-\end{prop}

\begin{Proof}
With \eqref{eq:maj2}, equation \eqref{eq:1}  becomes   
\begin{equation}  
\label{eq:D2}  
P(f(X)-E[f(X)]\geq x)\leq   
\exp\left(-\frac{(2-\alpha)R^\alpha}{2n\s }\right)+ 1-\exp\left(-\frac{\s}{\alpha R^\alpha}\right),  
\end{equation}  
as long as \eqref{eq:R-1} holds (still using \eqref{eq:Cond1} with  $y= R$).  
\vskip 5pt

We compare now the two summands in the right-hand side of \eqref{eq:D2}.  
To do so, let $u=\frac{(2-\alpha)R^\alpha}{2n\s }$ and let $c=\frac{2-\alpha}{2n\alpha}$ and compare $1-e^{-c/u}$ and $e^{-u}$. To this end, let us study for some $\varepsilon>0$ the function   
\begin{equation}  
\label{eq:h5}  
g_{n,\alpha}(u)=\varepsilon e^{-u}+e^{-c/u}  
\end{equation}  
and in particular let us see when it is larger than $1$. We cannot study the function $g_{n,\alpha}$ in \eqref{eq:h5} as easily as the function $g_{n,\alpha}$ in \eqref{eq:h4}, 
the argument is thus different from that of Section \ref{sec:intermediaire1}.  
However \eqref{eq:varep11} guarantees that $g_{n,\alpha}$ takes values larger than $1$
since  $g_{n,\alpha}(1)\geq c+e^{-c}\geq 1$, for all $c>0$.    
We investigate now for some interval containing $1$ where $g_{n,\alpha}$ is larger than $1$. Since $\varepsilon>ec$, we look for $c e^{1-u}+e^{-c/u}\geq 1$.  
   
Observe first that for some fixed $u$,  $c \to c e^{1-u}+e^{-c/u}\geq 1$  increases for $c \geq c_0(u)=u(u-1)-u\ln u>0$.  
  
Elementary computations show that $|u-1|^{3/2}\geq c_0(u)$ for all $u\leq u_1\simeq 3.2$ and   
$$  
|u-1|^{3/2}e^{1-u}+e^{-|u-1|^{3/2}/u}\geq 1  
$$  
for $0\leq u\leq u_2\simeq 1.68$. Next,   
for $\varepsilon >ec$, $|u-1|^{3/2}\leq c$, and $u\leq u_2$, we have   
$$  
1\leq |u-1|^{3/2}e^{1-u}+e^{-|u-1|^{3/2}/u}\leq c e^{1-u}+e^{-c/u}\leq \varepsilon e^{1-u}+e^{-c/u}.   
$$  
Finally, we derive   
$$  
1-\exp\left(-\frac{\s}{\alpha R^\alpha}\right)\leq \varepsilon \exp\left(-\frac{(2-\alpha)R^\alpha}{2n\s }\right),  
$$  
for $R$ in the range   
\begin{equation}  
\label{eq:rangeR2}  
0<\frac{2n\s}{2-\alpha}\left(1-\left(\frac {2-\alpha}{2n\alpha}\right)^{2/3}\right)  
<R^\alpha<  
\frac{2n\s}{2-\alpha}\left(1+\left(\left(\frac {2-\alpha}{2n\alpha}\right)^{2/3}\wedge 0.68\right)\right).  
\end{equation}  
We thus obtain a mobile range of varying length for which we have  
\begin{equation}  
\label{eq:22}  
P(f(X)-E[f(X)]\geq x)\leq (1+\varepsilon)\exp\left(-\frac{2-\alpha}{2n \s }R^\alpha\right).  
\end{equation}  
We now express the deviation \eqref{eq:22} and its conditions \eqref{eq:R-1} and \eqref{eq:rangeR2} in terms of $x$. Since   
\begin{equation*}  
\label{eq:R-x2}  
 R\leq x\leq  R\left(1+\frac {2-\alpha}{2n(\alpha-1)\left(1-\left(\frac {2-\alpha}{2n\alpha}\right)^{2/3}\right)}\right),  
\end{equation*}  
\eqref{eq:Deviation2} is obtained by using the bound \eqref{eq:22}.  
  
To express the range in terms of $x$, argue as in the proof of Proposition \ref{prop:D-s1.1} 
using the function $\theta$.   
Note that $\disp \s\leq \frac{2n\s}{2-\alpha}\left(1-\left(\frac {2-\alpha}{2n\alpha}\right)^{2/3}\right)$; indeed this is the same as   
$$  
\frac{2-\alpha}{2n}\leq 1-\left(\frac {2-\alpha}{2n\alpha}\right)^{2/3}  
$$  
and this is true since $\frac{2-\alpha}{2n\alpha} \leq 1/4$ and $\alpha u\leq 2u\leq 1-u^{2/3}$, for all $0 \leq u\leq 1/4$.   
  
We can thus use the function $\theta$ in the domain where it defines an increasing bijection, the conditions  \eqref{eq:R-1} and \eqref{eq:rangeR2} can be rewritten as:  
\begin{align*}  
&\theta\left(\frac{2n\s}{2-\alpha}\left(1-\left(\frac {2-\alpha}{2n\alpha}\right)^{2/3}\right)\right)\\  
&< x<  
 \theta\left(\frac{2n\s}{2-\alpha}\left(\left(1+\left(\left(\frac {2-\alpha}{2n\alpha}\right)^{2/3}\wedge 0.68\right)\right)\wedge u_n^*(\alpha)/2\right)\right).  
\end{align*}  
Finishing these computations, we obtain the range \eqref{eq:range-x2}.\CQFD
\end{Proof}


\begin{Rem}[Comments on the bounds in the domain \eqref{eq:range-x2}]
{\rm $ $  
  
We discuss below the behavior of the domain \eqref{eq:range-x2} as $\alpha$ goes to $2$ or as $n$ goes to $+\infty$. Let $\eta>0$.    Since $\disp \left(1+\frac {2-\alpha}{2n(\alpha-1)\left(1-\left(\frac {2-\alpha}{2n\alpha}\right)^{2/3}\right)}\right)^\alpha\leq 1+\eta$, the deviation \eqref{eq:Deviation2}  can be rewritten as  
$$  
P(f(X)-E[f(X)]\geq x)\leq (1+\varepsilon)\exp\left(-\frac {2-\alpha}{2(1+\eta)n\s} x^{\alpha}\right).  
$$   
Moreover, since $\left(1+\frac{2-\alpha}{2n(\alpha-1)\left( \log \frac {2n}{2-\alpha}\wedge u_n^*(\alpha)/2\right)}\right)^\alpha\geq 1$, the domain in \eqref{eq:range-x2} can thus be replaced, for $\alpha$ close enough to $2$ or for $n$ large enough, by   
\begin{align}  
\label{eq:range2.1}  
&\frac{2(1+\eta)n\s}{2-\alpha}\left(1-\left(\frac {2-\alpha}{2n\alpha}\right)^{2/3}\right)<x^\alpha< \\
\nonumber  
&\hskip 3cm 
\frac{2n\s}{2-\alpha}\left(\left(1+\left(\frac {2-\alpha}{2n\alpha}\right)^{2/3}\right)\wedge u_n^*(\alpha)/2\right).  
\end{align}  

\vskip 5pt  
  
But, we have seen that for $\alpha$ close enough to $2$, $u_n^*(\alpha)=u_n(\alpha)$   
 satisfies \eqref{eq:u2}, so that  we obtain 
 the following domain  with two different orders in the lower and upper bounds in $n$ or $\alpha$:
$$  
\frac{2(1+\eta)n\s}{2-\alpha}\left(1-\left(\frac {2-\alpha}{2n\alpha}\right)^{2/3}\right)< x^\alpha<  
\frac {2n \s }{2-\alpha} \log \sqrt\frac {n-1}{2-\alpha}.   
$$  
For $n$ large enough, $u_n^*(\alpha)\in[3,4\frac{3-\alpha}{2-\alpha}]$ is bounded so that the minimum in \eqref{eq:range2.1} is $u_n^*(\alpha)/2$
 and the domain is of the following type:   
$$  
\frac{2(1+\eta)n\s}{2-\alpha}\left(1-\left(\frac {2-\alpha}{2n\alpha}\right)^{2/3}\right)  
< x^\alpha<  \frac {3n \s }{2-\alpha}.   
$$  
Finally, observe that for $\alpha$ close enough to $2$ or for $n$ large enough,  $\varepsilon$ can be chosen arbitrary small in \eqref{eq:varep11}. We thus  recover the same type of deviation as that obtained in \cite{M2} for $\alpha$ close enough to $2$.   
  
}
\end{Rem}
\vskip 5pt

\begin{Rem}
{\rm $ $
\begin{itemize}
\item
The deviations \eqref{eq:Deviation1} and \eqref{eq:Deviation2} obtained in Sections \ref{sec:intermediaire1} and \ref{sec:intermediaire2} are of the same type. 
The ranges in \eqref{eq:range1-x}  and \eqref{eq:range-x2} cannot however be completely  compared since we cannot compare $u_1(n,\alpha)$ and $1-\left(\frac{2-\alpha}{2n\alpha}\right)^{2/3}$ ; 
this requires to study the sign of $h_{n,\alpha}\big(1-\left(\frac{2-\alpha}{2n\alpha}\right)^{2/3}\big)$.
In the limiting cases $\alpha\to 2$ or $n\to+\infty$, the ranges obtained is Section \ref{sec:intermediaire1} are better.

\item Three regimes of deviation for the Lipschitz function of stable vector $f(X)$ are thus available. 
Roughly speaking, the regimes and their ranges are the following:

\begin{itemize} 
\item for $x$ small, the order is $\exp\big(-cx^{\frac \alpha{\alpha-1}}\big)$;
\item for $x^\alpha$ of order $\s$, the order is $e^{-cx^\alpha}$;
\item for $x^\alpha$ of order bigger than $\s$, the order is $c/x^\alpha$.
\end{itemize}
Here $c$ stands for a generic constant, different at each occurence and which depends on $\alpha$ and the dimension $d$.

\end{itemize}
}
\end{Rem}
  

\section{Poisson space with a stable Lévy measure of index $\alpha\in(0,2)$}  
\label{sec:BHP}  
  
  In this part, we study the deviation $P(F-m(F)\geq x)$ where $F$ is a stochastic functional on the Poisson space $\Omega^{\rit^d}$ on $\rit^d$ equipped with the stable Lévy measure $\nu$ given by \eqref{eq:nu}.
  We recall that $\Omega^{\rit^d}$ denote the set of Radon measures   
$$
\Omega^{\rit^d} = \left\{   
 \omega = \sum_{i=1}^N \epsilon_{t_i} \ : \   
 (t_i)_{i=1}^{i=N} \subset \rit^d, \ t_i\not=t_j, \ \forall i\not= j, \   
 N\in \nit\cup \{ \infty \}\right\},
 $$ 
where $\epsilon_t$ denotes the Dirac measure at $t\in \rit^d$. In the sequel, $P$ is the Poisson measure with intensity $\nu$ on $\Omega^{\rit^d}$. On the Poisson space, we dispose of the linear, closable, finite difference operator   
$$   
D:L^2(\Omega^{\rit^d}, P)   \longrightarrow L^2(\Omega^{\rit^d}\times \rit^d,P\otimes \nu)  
$$   
defined via   
$$   
D_x F(\omega ) = F(\omega \cup \{ x \}) -F(\omega ), \ \ \   
 \ \ \ dP\times \nu (d\omega , dx)\mbox{-a.e.},   
$$   
 where as a convention we identify $\omega \in \Omega^{\rit^d}$ with its support,   
 cf. e.g. \cite{BHP} and the references therein.   
  
 In this section, we generalize the deviation bounds obtained in Section \ref{sec:M} to the case of stochastic functionals satisfying 
 $$  
|D_yF(\omega)|\leq   ||y||, \lng P(d\omega)\otimes \nu(dy)\mbox{-a.e.}  
$$
Roughly speaking, this hypothesis on $DF$ is the analogous of $f$ Lipschitz in Section \ref{sec:M}. Since the mean may not exist anymore, deviation are expressed henceforth with respect to a median $m(F)$ of the stochastic functional $F$. 
Note finally that the case of Lipschitz function of stable vector $f(X)$ can be recovered, as a particular case,  in this framework for any index $\alpha\in (0,2)$.

  As previously explained, configurations are truncated  and we deal with   
 the functional restricted to the truncated configuration $\omega_R$ with Lemma \ref{lemme:bis4} while the rest of the configuration $\omega_R^c$ is controlled by some function $\gamma$.   
\begin{eqnarray}  
\nonumber  
P(F-m(F)\geq x)&=&P(F-m(F)\geq x,\; \omega^c_R=\emptyset)+P(F-m(F)\geq x,\; \omega^c_R\not =\emptyset)\\  
\label{eq:tech2}  
& \leq &P(F_R-m(F)\geq x)+P(\{ \omega \in \Omega^X \ : \  \omega^c_R\not =\emptyset \} ).  
\end{eqnarray}  
The second summand in the right-hand side of \eqref{eq:tech2} is dominated as in (6.2) in \cite{BHP}:  
\begin{equation}  
\label{eq:gamma}  
P(\{ \omega \in \Omega^X \ : \    
 \omega\cap B_X (0,R)^c\not =\emptyset \})   
\leq \gamma (R)  
\end{equation}  
where $\gamma(R)$ can be chosed to be either $\gamma(R)=1-\exp\left(-\frac{\s}{\alpha R^\alpha}\right)$ or $\gamma(R)=\frac{\s}{\alpha R^\alpha}$. This choice will be discussed latter. 
Using the notations of the proof of Theorem 5.2 in \cite{BHP}, we have   
$$   
D_yg(F_R) (\omega ) \leq \vert D_y F(\omega_R ) \vert \leq    \vert y  
\vert_X, \quad P(d\omega)\otimes \nu(dy) \mbox{ a.e.}   
,   
$$   
 where $g(x) = (x-m(F_R))^+\wedge r$. Thus   
$$  
\sup_{y \in B_X(0,R)}  D_yg(F_R)\leq    R, \quad P\mbox{-a.s.}   
$$  
and for $k\geq 2$   
$$   
\|Dg(F_R)\|_{L^\infty(\Omega^X , L^k(\nu_R))}^k\leq   
 \frac{\s}{k-\alpha} R^{k-\alpha}.  
$$   
First, we have as in (5.6) of \cite{BHP}  
$$  
P(F_R - m(F_R)\geq y) \leq P(g(F_R)\geq y)\leq P(g(F_R)-E[g(F_R)]\geq y/2).  
$$  
We apply now the following lemma to $g(F_R)$, this lemma deals with functionals on truncated configurations $F_R(\omega)=F(\omega_R)$ and is the counterpart of Lemma \ref{lemme:bis3} in the same way Lemma 5.5 in \cite{BHP} is that of Lemma 2 in \cite{HM}.

\begin{lemma}   
\label{lemme:bis4}  
Let $R>0$, and let $F:\Omega^{\rit^d} \longrightarrow \rit$ such that:  
\begin{equation*}  
\label{eq:DyF}  
\sup_{y\in B(0,R)} \vert D_yF (\omega ) \vert\leq  R \quad P(d\omega )\mbox{-a.s.}  
\end{equation*}  
For any $n\geq 2$,  let $u_n(\alpha)$ be the unique non-zero (thus positive) solution of   
$$
\disp e^u-1-\frac{n-1}{2-\alpha} u=0,
$$   
and let   
$$  
u_n^*(\alpha)=\min_{1<k<\frac{n-1+\alpha}2}  \left(\left(k!\frac{(n-1)(k+1-\alpha)}{(n-k)(2-\alpha)}\right)^{1/(k-1)}\right)\wedge u_n(\alpha).  
$$  
Then for all   
\begin{equation}  
\label{eq:Cond2}  
x\leq x_0:= n\frac{\s R^{1-\alpha}}{2-\alpha} u_n^*(\alpha),  
\end{equation}   
 we have   
\begin{equation*}  
\label{eq:bis2-2}  
P(F_R-E[F_R] \geq x)\leq \exp\left(-\frac {(2-\alpha)x^2}{2n\s R^{2-\alpha}}\right).  
\end{equation*}  
\end{lemma}  
\begin{Proof}  
 The proof follows the lines of that of Lemma \ref{lemme:bis3} with the following changes. Starting from Proposition 2.2 in \cite{BHP} in place of Theorem 1 in \cite{H02} we have \eqref{eq:dev0}  with the function $h_R$ given by   
$$  
 h_R(s) = \sup_{(\om,\om')\in\Omega^X\times \Omega^X }\left| \int_{X}(e^{s D_yF(\om)}-1)\: D_{y}F(\om') \nu_R (dy)\right|, \quad s >0.
$$
Using the bounds $|D_yF |\leq K:= R$, for $y\in B(0,R)$ and for any $n \geq 2$ the inequality \eqref{eq:h-bound2}, we have   
\begin{align}   
\nonumber  
& h_R(s)\leq \sum_{k=1}^{n-1} \frac{s^k}{k!}\sup_{\omega, \omega'\in \Omega^X}   
\int_X |D_yF(\omega)|^k  |D_yF(\omega')|\; \nu_R(dy)\\  
\label {eq:alpha5}  
&\hskip 4cm +\frac{e^{sK}-\sum_{k=0}^{n-1}s^k K^k/k!}{K^n}\sup_{\omega, \omega'\in \Omega^X} \int_X  |D_yF(\omega)|^n |D_yF(\omega')| \; \nu_R(dy).   
\end{align}   
Thus, using the inequality $xy\leq x^p/p+y^q/q$, for $p^{-1}+q^{-1}=1$ and $x,y\geq 0$, we have  for $q=k+1$ and $p=\frac{k+1}k$, $1\leq k\leq n$:   
\begin{eqnarray*}  
\lefteqn{   
\sup_{\omega, \omega'\in \Omega^X} \int_X |D_yF(\omega)|^k  |D_yF(\omega')| \nu_R(dy)  
}   
\\  
& \leq & \frac k{k+1}\sup_{\omega, \omega'\in \Omega^X}   
\int_X |D_yF(\omega)|^{k+1}\nu_R(dy) +\frac 1{k+1}\sup_{\omega, \omega'\in \Omega^X}   
\int_X |D_yF(\omega')|^{k+1} \nu_R(dy)   
=\alpha_{k+1}  
,   
\end{eqnarray*}   
 where for all $k\geq 2$,   
\begin{eqnarray*}  
\alpha_k&=&\|DF\|_{L^\infty(\Omega^X, L^k(\nu_R))}^k   
\\   
 & \leq & \int_{\{|y|_2\leq R\}}|y|_2^k \nu_R(dy)  
\\   
 & = &  \int_{S^{d-1}} \sigma(d\xi) \int_{\{|r\xi|_2\leq R\}}r^{k-1-\alpha} dr\\  
&\leq& \frac{\s }{k-\alpha}R^{k-\alpha}:=\tilde \alpha_k.  
\end{eqnarray*}  
Plugging this last bound in \eqref{eq:alpha5}, we recover \eqref{eq:alpha4}  in the proof of Lemma \ref{lemme:bis3} and we finish similarly.
\CQFD  
\end{Proof}   
  
Applying Lemma \ref{lemme:bis4} fo $g(F_R)$, so that for $\disp   
y\leq 2x_0=2n \frac{\s  R^{1-\alpha}}{2-\alpha} u_n^*(\alpha)$  
we get   
\begin{eqnarray*}  
\nonumber  
P(F_R - m(F_R)\geq y) &\leq& P(g(F_R)-E[g(F_R)]\geq y/2) \\  
&\leq& \exp\left(-\frac{2-\alpha}{2n \s R^{2-\alpha}}y^2\right).   
\end{eqnarray*}  
Now, with $y= R$,   
\begin{equation}  
\label{eq:111}  
P(F_R - m(F_R)\geq  R) \leq \exp\left(-\frac{2-\alpha}{2n \s }R^\alpha \right)  
\end{equation}  
as long as, using  \eqref{eq:Cond2},  
\begin{equation*}  
\label{eq:R1}  
R^\alpha\leq  \frac {2n\s }{2-\alpha}  u_n^*(\alpha).  
\end{equation*}  
  
\vskip 5pt  
  
In order to control $P(F_R-m(F)\geq x)$ from (\ref{eq:111}),   
 we need to control $m(F_R)-m(F)$.  
To this end,  we apply Lemma 5.1 in \cite{BHP} with   
$$  
\beta(R) =2 R, \med   
\tilde{\gamma} (R) = \exp\left(-\frac{2-\alpha}{2n \s } R^{\alpha}\right),   
\med   
R_0= \left(\frac {2n\s }{2-\alpha}  u_n^*(\alpha)\right)^{1/\alpha}  
$$  
and where the hypothesis $R\geq R_0$  in Lemma 5.1 is replaced by $R\leq R_0$, so that the final condition on $R$ is changed in the same way.   
  
We thus have $m(F_R)-m(F)\leq  R$, for all $R$ such that   
\begin{equation}   
\label{ccond}   
\inf_{0<\delta < 1/2}   
 \max \left(   
 \gamma ^{-1}(\delta)   
 ,   
 \tilde{\gamma}^{-1}   
 \left( \frac{1}{2} - \delta \right)   
 \right)   
 \leq R\leq    
\left(\frac {2n\s }{2-\alpha}  u_n^*(\alpha)\right)^{1/\alpha},  
\end{equation}    
where $\gamma(R)$ is given in \eqref{eq:gamma}. We thus have for  $R=x/2$  
\begin{eqnarray*}  
P(F_R-m(F)\geq x)\leq P(F_R-m(F_R)\geq x/2)&\leq&   
\exp\left(-\frac{2-\alpha}{2n \s }R^\alpha \right)\\  
&\leq& \exp\left(-\frac{2-\alpha}{2n \s }\left(\frac x{2 }\right)^\alpha\right).  
\end{eqnarray*}  
 From \eqref{eq:tech2}, it follows that   
\begin{equation}  
\label{eq:D1}  
P(F_R-m(F)\geq x)\leq \exp\left(-\frac{2-\alpha}{2n \s }\left(\frac x{2 }\right)^\alpha\right)+ \gamma(x),
\end{equation}  
as long as the condition \eqref{ccond} holds, that is in terms of $x$ as long as   
\begin{equation}   
\label{ccondx}   
\inf_{0<\delta < 1/2}   
 \max \left(   
 \gamma ^{-1}(\delta)   
 ,   
 \tilde{\gamma}^{-1}   
 \left( \frac{1}{2} - \delta \right)   
 \right)   
 \leq \frac x{2 } \leq    
\left(\frac {2n\s }{2-\alpha} u_n^*(\alpha)\right)^{1/\alpha}.  
\end{equation}     
  
Next, we compare the two summands in the right-hand side of \eqref{eq:D1} using first $\gamma(R)=\frac{\s}{\alpha R^\alpha}$ and next $\gamma(R)=1-\exp\left(-\frac{\s}{\alpha R^\alpha}\right)$ in order to estimate the remainder term $P(\omega_R^c\not =0)$.   

  
\subsection*{First choice for the function $\gamma$}  
  
In this section, we take $\disp \gamma(R)=\frac{\s}{\alpha R^\alpha}$ and we have:
\begin{prop}  
\label{prop:D2-s}  
Let $F$ be a stochastic functional on the Poisson space equipped with the $\alpha$-stable Lévy measure \eqref{eq:nu}, $\alpha\in(0,2)$. Assume that   
$$  
|D_yF(\omega)|\leq   ||y||, \lng P(d\omega)\otimes \nu(dy)\mbox{-a.e.}  
$$  
Then for $\varepsilon$ satisfying \eqref{eq:varep1}, the following deviation holds 
\begin{equation}  
\label{eq:Dev2}  
P(F-m(F)\geq x)\leq (1+\varepsilon)\exp\left(-\frac{2-\alpha}{2n \s }\left(\frac x{2 }\right)^\alpha\right)  
\end{equation}  
for all $x$ in the range  
\begin{align}  
\label{eq:range2}   
&   
 \hskip -2cm    
 \left(\inf_{0<\delta < 1/2}   
 \max \left(   
 \gamma ^{-1}(\delta)   
 ,   
 \tilde{\gamma}^{-1}   
 \left( \frac{1}{2} - \delta \right)   
 \right)^{\alpha}\right)  
 \vee \frac {2n \s }{2-\alpha}u_1(n,\alpha)   
\\   
\nonumber   
& \hskip -0cm   
\leq \left(\frac x{2 }\right)^\alpha \leq  
  \frac {2n \s }{2-\alpha}\big(u_2(n,\alpha)\wedge u_n^*(\alpha)\big)
\end{align}    
where $u_1(n,\alpha)$, $u_1(n,\alpha)$ and $u_n^*(\alpha)$ are as in Proposition \ref{prop:D-s1.1} with $\delta=1$.
 \end{prop}    
This result generalizes Proposition 6.3 in \cite{BHP} in the same way Propositions \ref{prop:D-s1.1} and \ref{prop:D-s2.1} generalize Theorem 2 in \cite{M2}.
  
\begin{Proof}
As in Section \ref{sec:M}, we compare the two summands in \eqref{eq:D1} studying the function $h_{n,\alpha}$ in \eqref{eq:h4}.  
With notations as in Section \ref{sec:M} and for $\varepsilon$ satisfying \eqref{eq:varep1}, we derive \eqref{eq:Dev2} for
\begin{equation*}  
\label{eq:domaine-x1}  
\frac {2n \s }{2-\alpha}u_1(n,\alpha)< \left(\frac x{2 }\right)^\alpha<    
\frac {2n \s }{2-\alpha}u_2(n,\alpha),  
\end{equation*}  
and still under the condition \eqref{ccondx}.   \CQFD
\end{Proof}


\subsection*{Discussion on the range of deviation in \eqref{eq:range2}}

\begin{itemize}   
\item We have   
$$  
\gamma^{-1}(\delta)=\left(\frac{\s }{\alpha \delta}\right)^{1/\alpha},\med   
\tilde \gamma^{-1}(\delta)=\left(\frac{2n\s }{2-\alpha}\log (1/\delta)\right)^{1/\alpha}.  
$$  
Since there is a unique solution $\delta_0(n,\alpha)\in (0, 1/2)$ to the equation   
\begin{equation}  
\label{eq:az}  
\frac{2-\alpha}{2n \alpha}=\delta \log \frac 1{1/2-\delta},  
\end{equation}  
we have   
$\disp   
\inf_{0<\delta < 1/2}   
 \max \left(   
 \gamma ^{-1}(\delta)   
 ,   
 \tilde{\gamma}^{-1}   
 \left( \frac{1}{2} - \delta \right)   
 \right)^{\alpha}=\frac{\s }{\alpha\delta_0(n,\alpha)}$.  
  
\item Moreover, note that since the left-hand side of \eqref{eq:az} goes to zero when $\alpha\to 2$ or when $n\to +\infty$,  we have   
$\disp  
\delta_0(n,\alpha)\log \left(\frac 1{1/2-\delta_0(n,\alpha)}\right)\to 0  
$  
and thus $\delta_0(n,\alpha)\to 0$. It is also easy to deduce an equivalent for $\delta_0(n,\alpha)$ for $\alpha\to 2$ or for $n\to+\infty$:   
\begin{equation*}  
\delta_0(n,\alpha)\simeq \frac{2-\alpha}{2n\alpha\log 2}  
\end{equation*} 
where, when $x\to x_0$, $f(x)\simeq g(x)$ means $\disp \lim_{x\to x_0}\frac{f(x)}{g(x)}=1$.  
We thus have for $\alpha\to 2$ or for $n\to+\infty$,   
\begin{equation}  
\label{eq:lhs1}  
\inf_{0<\delta < 1/2}   
 \max \left(   
 \gamma ^{-1}(\delta)   
 ,   
 \tilde{\gamma}^{-1}   
 \left( \frac{1}{2} - \delta \right)   
 \right)^{\alpha}\simeq \frac {2n\s \log 2}{2-\alpha}.  
\end{equation}  
  
\item   
Since we have seen $\disp u_2(n,\alpha)>\log \left(\frac {2n\alpha \varepsilon}{2-\alpha}\right)>1$. The upper bound of the domain \eqref{eq:range2} in Proposition~\ref{prop:D2-s} can thus be taken to be   
\begin{equation}  
\label{eq:rhs}  
\frac {2n\s }{2-\alpha}  \left(\log \left(\frac {2n \alpha\varepsilon}{2-\alpha}\right) \wedge u_n^*(\alpha)\right).  
\end{equation}  
\end{itemize}


\subsection*{Comments on the bounds in the domain \eqref{eq:range2}}  
  
For $\alpha$ close enough to $2$ or for $n$ large enough, we have seen in Section \ref{sec:intermediaire1} the limits \eqref{eq:behavior}. Similarly, for any $\eta>0$ and for $\alpha$ close enough to $2$ or for $n$ large enough, the domain \eqref{eq:range2} reduces to   
\begin{equation*}  
(1+\eta)\s\left(\frac {2n \log 2}{2-\alpha}  
\vee \frac{\alpha}{\varepsilon}\right)< x^\alpha<  
\left\{  
\begin{array}{ll}  
\disp \frac {2n\s}{2-\alpha}\log\left(\sqrt\frac{n-1}{2-\alpha}\right)&\mbox{ if } \alpha\to 2,\\  
\disp \frac {6n\s}{2-\alpha} &\mbox{ if } n\to+\infty.  
\end{array}  
\right .  
\end{equation*}

\begin{Rem}  
{\rm  $ $ 
\begin{itemize}

\item
Proposition \ref{prop:D2-s} complements Theorem 6.1 in \cite{BHP} exhibiting for any $\alpha\in(0,2)$ a deviation regime in $\disp \exp\left(-\frac{2-\alpha}{2n\s }\left(\frac x{2 }\right)^\alpha\right)$  for $x$ in a finite range, while Theorem 6.1 gives a deviation regime in $1/x^\alpha$ for $x$ large enough, of order at least  
$$  
\s\frac 2{2-\alpha}\log \left(\frac 2{2-\alpha}\right)\log\left(\frac 2{2-\alpha}\log \left(\frac 2{2-\alpha}\right)\right).  
$$  

\item When $\alpha$ is close to $2$, we recover a range of the same type as in Theorem 6.3 in \cite{BHP}.

\end{itemize}
}\end{Rem}

  
\subsection*{Second choice for the function $\gamma$}  
  
Here, like in Section \ref{sec:intermediaire2}, we estimate the remainder term $P(\omega_R^c\not =0 )$ using directly \eqref{eq:maj2} rather than the bound \eqref{eq:maj3}, we obtain
the following result probably less sharp,  but with more explicit bounds. 
\begin{prop}  
\label{prop:D2-s2}  
Let $F$ be a stochastic functional on the Poisson space equipped with the $\alpha$-stable Lévy measure \eqref{eq:nu}, $\alpha\in(0,2)$. Assume that   
$$  
|D_yF(\omega)|\leq   ||y||, \lng P(d\omega)\otimes \nu(dy)\mbox{-a.e.}  
$$  
Then for any $\varepsilon$ satisfynig \eqref{eq:varep1}, the following  deviation holds 
\begin{equation}  
\label{eq:Dev3}  
P(F-m(F)\geq x)\leq (1+\varepsilon)\exp\left(-\frac{2-\alpha}{2n \s }\left(\frac x{2 }\right)^\alpha\right)  
\end{equation}  
for all $x$ in the range  
\begin{align}  
\label{eq:range23}   
&   
 \hskip -2cm    
 \left(\inf_{0<\delta < 1/2}   
 \max \left(   
 \gamma ^{-1}(\delta)   
 ,   
 \tilde{\gamma}^{-1}   
 \left( \frac{1}{2} - \delta \right)   
 \right)^{\alpha}\right)  
 \vee \frac{2n\s}{2-\alpha}\left(1-\left(\frac {2-\alpha}{2n\alpha}\right)^{2/3}\right)  
\\   
\nonumber   
& \hskip -0cm  <\left(\frac{x}{2 }\right)^\alpha<  
\frac{2n\s}{2-\alpha}\left(\left(1+\left(\left(\frac {2-\alpha}{2n\alpha}\right)^{2/3}\wedge 0.68\right)\right)\wedge u_n^*(\alpha)/2\right).  
\end{align}    
\end{prop}     
Taking limits in \eqref{eq:range23} yields nothing since $1\pm(\frac{2-\alpha}{2n\alpha})^{2/3} \to 0$ when $\alpha \to 2$ or when $n\to+\infty$.   

\begin{Proof}  
As in Section \ref{sec:intermediaire2}, we compare the two summands in \eqref{eq:D1} studying now the function $g_{n,\alpha}$ in \eqref{eq:h5}.  
As in Section \ref{sec:intermediaire2}, we derive   \eqref{eq:Dev3} for $x$ in the range  
\begin{equation*}  
\label{eq:rangeR3}  
0<\frac{2n\s}{2-\alpha}\left(1-\left(\frac {2-\alpha}{2n\alpha}\right)^{2/3}\right)  
<\left(\frac x{2 }\right)^\alpha<  
\frac{2n\s}{2-\alpha}\left(1+\left(\left(\frac {2-\alpha}{2n\alpha}\right)^{2/3}\wedge 0.68\right)\right),  
\end{equation*}  
and still under condition \eqref{ccondx}.\CQFD  
\end{Proof}

  
\subsection*{Lower range for index $\alpha>1$}

Finally, note that it is a simple verification that the counterpart of Proposition~\ref{prop:dev-petit2}  holds for functional on
the Poisson space, namely,
\begin{prop}  
Let $F$ be a stochastic functional on the Poisson space equipped with the stable Lévy measure \eqref{eq:nu} with index $\alpha\in(1,2)$. Assume that   
$$  
|D_yF(\omega)|\leq   ||y||, \lng P(d\omega)\otimes \nu(dy)\mbox{-a.e.}  
$$  
Then for all $n$ large enough and for all
$$
\lambda \in \left]\frac{\s}{\alpha-1}\left(1+\sqrt{\frac{2(\alpha-1)^2n \s^{\frac{2-\alpha}{\alpha-1}}}{\alpha(2-\alpha)}}\right), \frac \s{\alpha-1}\left(1+\frac{\alpha-1}{2-\alpha}n u_n^*(\alpha)\right)\right[,
$$
there exists $x_1(n,\alpha, \lambda)>0$ such that for all $0\leq x\leq  x_1(n,\alpha,\lambda)$, 
\begin{align*}
\nonumber
&P(f(X)-E[f(X)]\geq x)\\
&\leq 
 \exp\left(-\frac{2-\alpha}{2n\s^{1/(\alpha-1)}}(\lambda-\frac \s{\alpha-1})^2
\left(\frac x\lambda \right)^{\frac \alpha{\alpha-1}} \right)
+\frac \s\alpha \left(\frac x\lambda \right)^{\frac \alpha{\alpha-1}}\leq 1.
\end{align*}
\end{prop}  

The same comments as in Remarks \ref{rem:petit1} and \ref{rem:petit3} apply. 

For stable index $\alpha\leq 1$, lower range deviation requires to investigate the convergence of the median of truncated functional which, at this point, requires more work.




\begin{thebibliography}{9}  
  
\bibitem{BHP} J.-C. {\sc Breton}, C. {\sc Houdr\'e},  N. {\sc Privault}, {\it Dimension free and infinite variance tail estimates on Poisson space}, Preprint, 2004.   
  
\bibitem{H02} C. {\sc Houdr\'e}, {\it Remarks on deviation inequalities for functions of infinitely divisible random vectors}, 
Ann. Probab., vol. 30, no. 3, pp. 1123--1237, 2002.

\bibitem{HM} C. {\sc Houdr\'e}, P. {\sc Marchal}, {\it On the concentration of measure phenomenon for stable and related random vectors}, Ann. Probab., vol. 32, no. 2, p. 1496--1508, 2004.  
   
\bibitem{M2} P. {\sc Marchal}, {\it A note on measure concentration for stable distributions with index close to $2$}, Preprint, 2004.   
   
\end{thebibliography}
\end{document}